\documentclass[reqno]{amsart}
\usepackage{amssymb}
\usepackage[colorlinks=true,hypertex]{hyperref}
\allowdisplaybreaks[4]
\theoremstyle{plain}
\newtheorem{thm}{Theorem}

\theoremstyle{remark}
\newtheorem{rem}{Remark}

\date{Commenced on 29 January 2009 and completed on 16 February 2009 in Melbourne}
\date{}

\begin{document}

\title[Concise sharpening and generalizations of Shafer's inequality]
{Concise sharpening and generalizations of Shafer's inequality for the arc sine function}

\author[F. Qi]{Feng Qi}
\address[F. Qi]{Research Institute of Mathematical Inequality Theory, Henan Polytechnic University, Jiaozuo City, Henan Province, 454010, China}
\email{\href{mailto: F. Qi <qifeng618@gmail.com>}{qifeng618@gmail.com}, \href{mailto: F. Qi <qifeng618@hotmail.com>}{qifeng618@hotmail.com}, \href{mailto: F. Qi <qifeng618@qq.com>}{qifeng618@qq.com}}
\urladdr{\url{http://qifeng618.spaces.live.com}}

\author[B.-N. Guo]{Bai-Ni Guo}
\address[B.-N. Guo]{School of Mathematics and Informatics, Henan Polytechnic University, Jiaozuo City, Henan Province, 454010, China}
\email{\href{mailto: B.-N. Guo <bai.ni.guo@gmail.com>}{bai.ni.guo@gmail.com}, \href{mailto: B.-N. Guo <bai.ni.guo@hotmail.com>}{bai.ni.guo@hotmail.com}}
\urladdr{\url{http://guobaini.spaces.live.com}}

\begin{abstract}
In this paper, by a concise and elementary approach, we sharpen and generalize Shafer's inequality for the arc sine function, and some known results are extended and generalized.
\end{abstract}

\keywords{sharpening, generalization, Shafer's inequality, arc sine function, monotonicity}

\subjclass[2000]{Primary 33B10; Secondary 26D05}

\thanks{The first author was partially supported by the China Scholarship Council}

\thanks{This paper was typeset using \AmS-\LaTeX}

\maketitle

\section{Introduction and main results}

In~\cite[p.~247, 3.4.31]{mit}, it was listed that the inequality
\begin{equation}\label{Shafer-ineq-arcsin}
\arcsin x>\frac{6\bigl(\sqrt{1+x}\,-\sqrt{1-x}\,\bigr)}{4+\sqrt{1+x}\,+\sqrt{1-x}\,} >\frac{3x}{2+\sqrt{1-x^2}\,}
\end{equation}
holds for $0<x<1$. It was also pointed out in~\cite[p.~247, 3.4.31]{mit} that these inequalities are due to R. E. Shafer, but no a related reference is cited. By now we do not know the very original source of inequalities in~\eqref{Shafer-ineq-arcsin}.
\par
In~\cite{Zhu-Mia-05-Shafer}, the left-hand side inequality in~\eqref{Shafer-ineq-arcsin} was recovered and an upper bound was presented as follows:
\begin{equation}\label{Zhu-Mia-05-Shafer-inq}
\arcsin x \le \frac{\pi\bigl(\sqrt2\,+1/2\bigr)\bigl(\sqrt{1+x}\,-\sqrt{1-x}\,\bigr)} {4+\sqrt{1+x}\,+\sqrt{1-x}\,}, \quad 0\le x\le1.
\end{equation}
\par
In \cite{Malesevic-JIA-07-78691, Malesevic-MIA-07-3, Zhu-Shafer-JIA-07-67430}, the upper bound in~\eqref{Zhu-Mia-05-Shafer-inq} was numerically improved to
\begin{equation}\label{Malesevic-ineq-improve}
\arcsin x\le \frac{\bigl[\pi\bigl(2-\sqrt2\,\bigr)\big/\bigl(\pi-2\sqrt2\,\bigr)\bigr] \bigl(\sqrt{1+x}\,-\sqrt{1-x}\,\bigr)} {(4-\pi)\sqrt2\,\big/\bigl(\pi-2\sqrt2\,\bigr)+\sqrt{1+x}\,+\sqrt{1-x}\,},\quad 0\le x\le1.
\end{equation}
\par
For more information, please refer to~\cite{Oppeheim-Sin-Cos.tex} and related references therein.
\par
The aim of this paper is to sharpen and generalize the above inequalities.
\par
Our main results may be stated as follows.

\begin{thm}\label{Shafer-ArcSin-thm}
Let $\alpha$ be a real number. Then the function
\begin{equation}
f_\alpha(x)=\frac{\alpha+\sqrt{1+x}\,+\sqrt{1-x}\,}{\sqrt{1+x}\,-\sqrt{1-x}\,}\arcsin x, \quad x\in(0,1]
\end{equation}
has the following properties:
\begin{enumerate}
  \item
For $\alpha\ge4$, it is strictly increasing;
  \item
For $\alpha\le\frac{4(\pi-2)}{\sqrt{2}\,(4-\pi)}$, it is strictly decreasing;
\item
For $4>\alpha>\frac{4(\pi-2)}{\sqrt{2}\,(4-\pi)}$, it has a unique minimum.
\end{enumerate}
\end{thm}

As direct consequences of Theorem~\ref{Shafer-ArcSin-thm}, the following inequalities may be derived.

\begin{thm}\label{Shafer-ArcSin-thm-2}
If $\alpha\ge4$, then the inequality
\begin{multline}\label{Shafer-ArcSin-thm-2-ineq}
\frac{(2+\alpha)\bigl(\sqrt{1+x}\,-\sqrt{1-x}\,\bigr)}{\alpha+\sqrt{1+x}\,+\sqrt{1-x}\,}<\arcsin x \\*
<\frac{\bigl[\pi\bigl(\sqrt2\,+\alpha\bigr)/2\sqrt2\,\bigr] \bigl(\sqrt{1+x}\,-\sqrt{1-x}\,\bigr)} {\alpha+\sqrt{1+x}\,+\sqrt{1-x}\,},\quad x\in(0,1).
\end{multline}
If $4>\alpha>\frac{4(\pi-2)}{\sqrt{2}\,(4-\pi)}$, then
\begin{equation}\label{Shafer-ArcSin-thm-2-ineq-2}
\arcsin x <\frac{\max\bigl\{2+\alpha,\pi\bigl(\sqrt2\,+\alpha\bigr)/2\sqrt2\,\bigr\} \bigl(\sqrt{1+x}\,-\sqrt{1-x}\,\bigr)} {\alpha+\sqrt{1+x}\,+\sqrt{1-x}\,},\quad x\in(0,1).
\end{equation}
If $\alpha\le\frac{4(\pi-2)}{\sqrt{2}\,(4-\pi)}$, then the
inequality~\eqref{Shafer-ArcSin-thm-2-ineq} reverses.
\par
Moreover, the constants $2+\alpha$ and $\frac{\pi(\sqrt2\,+\alpha)}{2\sqrt2\,}$ in~\eqref{Shafer-ArcSin-thm-2-ineq} and the scalar $\max\Bigl\{2+\alpha,\frac{\pi(\sqrt2\,+\alpha)}{2\sqrt2\,}\Bigr\}$ in~\eqref{Shafer-ArcSin-thm-2-ineq-2} are the best possible.
\end{thm}

\begin{rem}
It is easy to see that the left-hand side inequality in~\eqref{Shafer-ineq-arcsin} can be deduced from the left-hand inequality in~\eqref{Shafer-ArcSin-thm-2-ineq} by taking $\alpha=4$, that the inequality~\eqref{Zhu-Mia-05-Shafer-inq} is the special case $\alpha=4$ of the right-hand side inequality in~\eqref{Shafer-ArcSin-thm-2-ineq}, and that the inequality~\eqref{Malesevic-ineq-improve} is the special case $\alpha=\frac{(4-\pi)\sqrt2\,}{\pi-2\sqrt2\,}$ of the inequality in~\eqref{Shafer-ArcSin-thm-2-ineq-2}. Therefore, our Theorem~\ref{Shafer-ArcSin-thm} and Theorem~\ref{Shafer-ArcSin-thm-2} extend, sharpen and generalize related results demonstrated in~\cite{Malesevic-JIA-07-78691, Malesevic-MIA-07-3, mit, Zhu-Mia-05-Shafer, Zhu-Shafer-JIA-07-67430}.
\end{rem}

\begin{rem}
Comparing with the methods used in~\cite{Malesevic-JIA-07-78691, Malesevic-MIA-07-3, Zhu-Mia-05-Shafer, Zhu-Shafer-JIA-07-67430}, not only our proofs for Theorem~\ref{Shafer-ArcSin-thm} and Theorem~\ref{Shafer-ArcSin-thm-2} are more elementary and concise, but also we procure more general conclusions.
\end{rem}

\begin{rem}
For $4>\alpha>\frac{4(\pi-2)}{\sqrt{2}\,(4-\pi)}$, can one give a lower bound of the inequality~\eqref{Shafer-ArcSin-thm-2-ineq-2}?
\end{rem}

\section{Proofs of theorems}

Now we are in a position to prove our theorems.

\begin{proof}[Proof of Theorem~\ref{Shafer-ArcSin-thm}]
For $x\in(0,1)$, direct differentiation yields
\begin{gather*}
f_\alpha'(x)=\frac{\bigl[\alpha\bigl(\sqrt{1-x}\,+\sqrt{x+1}\,\bigr)+4\bigr]\sqrt{1-x^2}\,} {4\bigl(1-x^2\bigr) \bigl(1-\sqrt{1-x^2}\,\bigr)}\\
\times\Biggl\{\frac{2 \bigl\{2x\sqrt{1-x^2}\,+\alpha \bigl[x\bigl(\sqrt{1-x}\,+\sqrt{x+1}\,\bigr) +\sqrt{1-x}\,-\sqrt{x+1}\,\bigr] \bigr\}} {\bigl[\alpha \bigl(\sqrt{1-x}\,+\sqrt{x+1}\,\bigr)+4\bigr]\sqrt{1-x^2}\,}- \arcsin x\Biggr\}\\
\triangleq\frac{\bigl[\alpha\bigl(\sqrt{1-x}\,+\sqrt{x+1}\,\bigr)+4\bigr]\sqrt{1-x^2}\,} {4\bigl(x^2-1\bigr) \bigl(\sqrt{1-x^2}\,-1\bigr)}h_\alpha(x),\\
h_\alpha'(x)=\frac{2\bigl(x^2+\sqrt{1-x^2}\,-1\bigr)}{\bigl(1-x^2\bigr)\bigl[\alpha \bigl(\sqrt{1-x}\,+\sqrt{x+1}\,\bigr)+4\bigr]^2} \biggl\{\alpha^2-8\\
+\frac{x\bigl(\sqrt{1-x}\,-\sqrt{x+1}\,\bigr)\bigl(\sqrt{1-x^2}\,-2\bigr) +2\bigl(\sqrt{1-x}\,+\sqrt{x+1}\,\bigr)   \bigl(\sqrt{1-x^2}\,-1\bigr)}{x^2+\sqrt{1-x^2}\,-1} \alpha\biggr\}\\
\triangleq\frac{2\bigl(x^2+\sqrt{1-x^2}\,-1\bigr)\bigl[\alpha^2-8+\alpha g(x)\bigr]} {\bigl(1-x^2\bigr)\bigl[\alpha \bigl(\sqrt{1-x}\,+\sqrt{x+1}\,\bigr)+4\bigr]^2},\\
\begin{aligned}
g'(x)&=-\frac{2 \bigl(\sqrt{1-x}\,-\sqrt{x+1}\,\bigr)+x\bigl(\sqrt{1-x}\,+\sqrt{x+1}\,\bigr)} {2\bigl(x^2+\sqrt{1-x^2}\,-1\bigr)}\\
&\triangleq-\frac{p(x)}{2\bigl(x^2+\sqrt{1-x^2}\,-1\bigr)},\\
p'(x)&=\frac{3x\bigl(\sqrt{1-x}\,-\sqrt{x+1}\,\bigr)}{2\sqrt{1-x^2}\,}\\
&<0.
\end{aligned}
\end{gather*}
Since $p(0)=0$ and $p(x)$ is strictly decreasing on $(0,1)$, the function $p(x)$ is negative on $(0,1)$, so the derivative $g'(x)$ is positive and $g(x)$ is strictly increasing on $(0,1)$.
\par
By virtue of
$$
\lim_{x\to0^+}g(x)=-2\quad\text{and}\quad \lim_{x\to1^-}g(x)=-\sqrt2\,,
$$
it follows for $\alpha>0$ that
\begin{enumerate}
\item
when $\alpha^2-2\alpha-8\ge0$, that is, $\alpha\ge4$, the derivative $h_\alpha'(x)$ is positive and the function $h_\alpha(x)$ is strictly increasing on $(0,1)$;
\item
when $\alpha^2-\sqrt2\,\alpha-8\le0$, that is, $0<\alpha\le\frac{\sqrt2\,+\sqrt{34}\,}2$, the derivative $h_\alpha'(x)$ is negative and the function $h_\alpha(x)$ is strictly decreasing on $(0,1)$;
\item
when $\alpha^2-2\alpha-8<0$ and $\alpha^2-\sqrt2\,\alpha-8>0$, that is, $4>\alpha>\frac{\sqrt2\,+\sqrt{34}\,}2$, the derivative $h_\alpha'(x)$ has a unique zero and the function $h_\alpha(x)$ has a unique minimum on $(0,1)$.
\end{enumerate}
It is easy to see that
\begin{equation}\label{h-alpha-limits}
\lim_{x\to0^+}h_\alpha(x)=0\quad\text{and}\quad \lim_{x\to1^-}h_\alpha(x)=\frac{8-4\pi-\sqrt{2}\,(\pi-4)\alpha}{2\sqrt{2}\,\alpha+8}.
\end{equation}
Whence the function $h_\alpha(x)$ and $f_\alpha'(x)$ are strictly positive for
$\alpha\ge4$ and strictly negative for $0<\alpha\le\frac{4(\pi-2)}{\sqrt{2}\,(4-\pi)}$
on $(0,1)$ and have a unique zero on $(0,1)$ for
$4>\alpha>\frac{4(\pi-2)}{\sqrt{2}\,(4-\pi)}$. Consequently, the function $f_\alpha(x)$
is strictly increasing for $\alpha\ge4$, strictly decreasing for
$0<\alpha\le\frac{4(\pi-2)}{\sqrt{2}\,(4-\pi)}$, and has a unique minimum for
$4>\alpha>\frac{4(\pi-2)}{\sqrt{2}\,(4-\pi)}$.
\par
For $x\in(0,1)$ and $\alpha<0$, considering the fact that the function $\alpha
\bigl(\sqrt{1-x}\,+\sqrt{x+1}\,\bigr)+4$ for $x\in(0,1)$ does not equal zero if and
only if $0>\alpha\ge-2$ or $\alpha\le-2\sqrt2\,=-2.828\dotsm$, similar argument as
above gives that
\begin{enumerate}
\item
when $\alpha^2-2\alpha-8\le0$, that is, $0>\alpha\ge-2$, the derivative $h_\alpha'(x)$ is negative and the function $h_\alpha(x)$ is strictly decreasing on $(0,1)$;
\item
when $\alpha^2-\sqrt2\,\alpha-8\ge0$, that is, $\alpha\le-2\sqrt2\,<\frac{\sqrt2\,-\sqrt{34}\,}2=-2.208\dotsm$, the derivative $h_\alpha'(x)$ is positive and the function $h_\alpha(x)$ is strictly increasing on $(0,1)$;
\item
when $\alpha^2-2\alpha-8>0$ and $\alpha^2-\sqrt2\,\alpha-8<0$, that is, $-2>\alpha>\frac{\sqrt2\,-\sqrt{34}\,}2>-2\sqrt2\,$, the derivative $h_\alpha'(x)$ has a unique zero and the function $h_\alpha(x)$ has a unique maximum on $(0,1)$, and so, by~\eqref{h-alpha-limits}, the function $h_\alpha(x)$ has a unique zero on $(0,1)$.
\end{enumerate}
As a result, from the first limit in~\eqref{h-alpha-limits}, it follows that the
function $h_\alpha(x)$ is negative for $0>\alpha\ge-2$ and positive for
$\alpha\le-2\sqrt2\,$. So $f'(x)$ is strictly negative for $0>\alpha\ge-2$ or
$\alpha\le-2\sqrt2\,$ on $(0,1)$. In a word, the function $f(x)$ is strictly decreasing
for $0>\alpha\ge-2$ or $\alpha\le-2\sqrt2\,$.
\par
For $x\in(0,1)$ and $\alpha=0$, the derivative of $f_\alpha(x)$ equals
$$
f_\alpha'(x)=\frac{x-\arcsin x}{x^2-1+\sqrt{1-x^2}}<0,\quad x\in(0,1).
$$
Thus, the function $f_0(x)$ is strictly decreasing on $(0,1)$.
\par
On the other hand, the derivative $f_\alpha'(x)$ may be rewritten as
\begin{align*}
f_\alpha'(x)&=\frac{\sqrt{1-x^2}\,} {4\bigl(1-x^2\bigr) \bigl(1-\sqrt{1-x^2}\,\bigr)}\\
&\quad\times\Biggl\{\frac{2 \bigl\{2x\sqrt{1-x^2}\, +\alpha
\bigl[x\bigl(\sqrt{1-x}\,+\sqrt{x+1}\,\bigr) +\sqrt{1-x}\,-\sqrt{x+1}\,\bigr] \bigr\}}
{\sqrt{1-x^2}\,}\\
&\quad-\bigl[\alpha\bigl(\sqrt{1-x}\,+\sqrt{x+1}\,\bigr)+4\bigr] \arcsin x\Biggr\}\\
&\triangleq \frac{\sqrt{1-x^2}\,} {4\bigl(1-x^2\bigr)
\bigl(1-\sqrt{1-x^2}\,\bigr)}F_\alpha(x)
\end{align*}
with
\begin{equation*}
F_\alpha'(x)=\frac{\alpha\sqrt{1-x^2}\,\bigl(\sqrt{1-x}\,-\sqrt{x+1}\,\bigr) \arcsin
x+8 \bigl(x^2-1+\sqrt{1-x^2}\,\bigr)}{2\bigl(x^2-1\bigr)}.
\end{equation*}
It is clear that when $\alpha\le0$ the derivative $F_\alpha'(x)$ is negative on
$(0,1)$, and so the function $F_\alpha(x)$ is strictly decreasing. By virtue of
$\lim_{x\to0^+}F_\alpha(x)=0$, it is deduced that $F_\alpha(x)<0$ on $(0,1)$, which
means that the function $f_\alpha(x)$ is strictly decreasing on $(0,1)$. The proof of
Theorem~\ref{Shafer-ArcSin-thm} is complete.
\end{proof}

\begin{proof}[Proof of Theorem~\ref{Shafer-ArcSin-thm-2}]
It is easy to obtain that
$$
\lim_{x\to0^+}f_\alpha(x)=2+\alpha\quad\text{and} \quad
\lim_{x\to1^-}f_\alpha(x)=\frac{\pi\bigl(\sqrt2\,+\alpha\bigr)}{2\sqrt2\,}.
$$
Hence, when $\alpha\ge4$, it follows from the increasing monotonicity in
Theorem~\ref{Shafer-ArcSin-thm} that
$$
2+\alpha<f_\alpha(x)<\frac{\pi\bigl(\sqrt2\,+\alpha\bigr)}{2\sqrt2\,}
$$
which can be rearranged as the inequality~\eqref{Shafer-ArcSin-thm-2-ineq}.
\par
When $\alpha\le\frac{4(\pi-2)}{\sqrt{2}\,(4-\pi)}$, the reversed version
of~\eqref{Shafer-ArcSin-thm-2-ineq} follows easily from the decreasing monotonicity of
the function $f_\alpha(x)$ presented in Theorem~\ref{Shafer-ArcSin-thm}.
\par
By the proof of Theorem~\ref{Shafer-ArcSin-thm}, when $4>\alpha>\frac{4(\pi-2)}{\sqrt{2}\,(4-\pi)}$, the function $f_\alpha(x)$ has a unique minimum on $(0,1)$, which means that
$$
f_\alpha(x)<\max\Bigl\{\lim_{x\to0^+}f_\alpha(x), \lim_{x\to1^-}f_\alpha(x)\Bigr\},\quad x\in(0,1).
$$
Rearranging this inequality yields~\eqref{Shafer-ArcSin-thm-2-ineq-2}.
The proof of Theorem~\ref{Shafer-ArcSin-thm-2} follows.
\end{proof}

\begin{rem}
This paper is a slightly modified version of the
preprint~\cite{Shafer-ArcSin.tex-arXiv}.
\end{rem}


\begin{thebibliography}{9}

\bibitem{Malesevic-JIA-07-78691}
B. J. Male\v{s}evi\'c, \emph{One method for proving inequalities by computer}, J. Inequal. Appl. \textbf{2007} (2007), Article ID 78691, 8~pages; Available online at \url{http://dx.doi.org/10.1155/2007/78691}.

\bibitem{Malesevic-MIA-07-3}
B. J. Male\v{s}evi\'c, \emph{An application of $\lambda$-method on inequalities of Shafer-Fink's type}, Math. Inequal. Appl. \textbf{10} (2007), no.~3, 529\nobreakdash--534.

\bibitem{mit}
D. S. Mitrinovi\'c, \textit{Analytic Inequalities}, Springer-Verlag, 1970.

\bibitem{Oppeheim-Sin-Cos.tex}
F. Qi and B.-N. Guo, \emph{A concise proof of Oppenheim's double inequality relating to the cosine and sine functions}, Avaliable online at \url{http://arxiv.org/abs/0902.2511}.

\bibitem{Shafer-ArcSin.tex-arXiv}
F. Qi and B.-N. Guo, \textit{Concise sharpening and generalizations of Shafer's
inequality for the arc sine function}, Available online at
\url{http://arxiv.org/abs/0902.2588}.

\bibitem{Zhu-Mia-05-Shafer}
L. Zhu, \emph{On Shafer-Fink inequalities}, Math. Inequal. Appl. \textbf{8} (2005), no.~4, 571\nobreakdash--574.

\bibitem{Zhu-Shafer-JIA-07-67430}
L. Zhu, \emph{On Shafer-Fink-type inequality}, J. Inequal. Appl. \textbf{2007} (2007), Article ID 67430, 4~pages; Avaliable online at \url{http://dx.doi.org/10.1155/2007/67430}.

\end{thebibliography}
\end{document}